\documentclass[lettersize,journal]{article}

\usepackage{arxiv}

\usepackage{amssymb,amsmath,amsthm,amsthm}
\allowdisplaybreaks 
\usepackage{amsbsy}
\usepackage{mathtools}
\usepackage{enumitem}
\usepackage{algorithm}
\usepackage{array}
\usepackage{textcomp}
\usepackage{stfloats}
\usepackage{url}
\usepackage{verbatim}
\usepackage{graphicx}
\usepackage{cite}
\usepackage{multirow}
\usepackage{tikz}
\usepackage{soul}
\usepackage[hypertexnames=false,hidelinks]{hyperref}
\usepackage{subcaption}
\usepackage{caption}
\usepackage{todonotes}
\usepackage{booktabs}
\usepackage{multirow, makecell}
\usepackage{xurl}

\newtheoremstyle{nonitalic} 
{3pt} 
{3pt} 
{\upshape} 
{} 
{\bfseries} 
{.} 
{.5em} 
{} 

\theoremstyle{nonitalic}
\theoremstyle{nonitalic}

\theoremstyle{plain} 

\usepackage{booktabs}

\newcommand{\R}{\mathbb{R}}

\newcommand{\prox}{\mathrm{prox}}

\newcommand{\be}{{\boldsymbol{e}}}
\newcommand{\bw}{{\boldsymbol{w}}}

\newcommand{\bx}{{\boldsymbol{x}}}
\newcommand{\by}{{\boldsymbol{y}}}
\newcommand{\bz}{{\boldsymbol{z}}}

\newcommand{\TpV}{T{\it p}V}

\usepackage{xcolor}

\usepackage{tcolorbox}
\newtcolorbox{low-dose-CA-problem}{
	colback=azure!5!white,
	colframe=azure!75!black,
	title=The Low-to-High Contrast Agent  Problem
}

\usepackage{tikz}
\usepackage{tikz-network}
\usetikzlibrary{3d,arrows.meta}

\date{}

\begin{document}
\title{Deep Guess acceleration for explainable image reconstruction in sparse-view CT}

\author{
Elena Loli Piccolomin \\
Department of Computer Science and Engineering, \\ University of Bologna,\\ Bologna, 40126, Italy. \\ \texttt{elena.loli@unibo.it}. \\
\And Davide Evangelista \\
Department of Computer Science and Engineering, \\ University of Bologna,\\ Bologna, 40126, Italy. \\
\texttt{davide.evangelista5@unibo.it} \\
\And Elena Morotti \\
Department of Political and Social Sciences, \\ University of Bologna,\\ Bologna, 40126, Italy. \\ \texttt{elena.morotti4@unibo.it}. 
}

\maketitle

\begin{abstract}
Sparse-view Computed Tomography (CT) is an emerging protocol designed to reduce X-ray dose radiation in medical imaging. Traditional Filtered Back Projection  algorithm reconstructions suffer from severe artifacts due to sparse data. In contrast, Model-Based Iterative Reconstruction (MBIR) algorithms, though better at mitigating noise through regularization, are too computationally costly for clinical use. 
This paper introduces a novel technique, denoted as the Deep Guess acceleration scheme, using a trained neural network both to quicken the regularized MBIR and to enhance the reconstruction accuracy.
We integrate state-of-the-art deep learning tools to initialize a clever starting guess for a proximal algorithm solving a non-convex model and thus computing an interpretable solution image in a few iterations.
Experimental results on real CT images demonstrate the Deep Guess effectiveness in (very) sparse tomographic protocols, where it overcomes its mere variational counterpart and many data-driven approaches at the state of the art.
We also consider a ground truth-free implementation and test the robustness of the proposed framework to noise.
\end{abstract}


\section{Introduction}\label{sec:introduction}

Computed tomography (CT) reconstruction from sparse views has emerged as a focal point in medical imaging research. A key challenge lies in diminishing the total X-ray radiation dose administered to patients while maintaining reconstruction precision. The sparse-view methodology revolves around decreasing the number of rotation angles utilized, thereby eliciting compromised image quality due to the introduction of numerous streaking artifacts by the classical analytical algorithms, such as the Filtered Back Projection (FBP).

Model-Based Iterative Reconstruction (MBIR) methods are extensively employed, relying on mathematical models to represent the problem being solved. In this case, sparse CT image reconstruction can be addressed as a linear problem of the form:
\begin{equation}
    \by = A\bx+ \be^\nu, 
    \label{eq:syslin}
\end{equation}
where  $\by\in \R^m$ is the sinogram, affected by noise  $\be^\nu \in \R^m$ with intensity depending on the positive parameter $\nu$, and $A \in \R^{m \times n}$ is the forward operator describing the CT system geometry. Due to the sparse sampling, $m<n$ and the system \eqref{eq:syslin} has infinitely many solutions. Additionally, owing to the typical ill-posedness of inverse problems, it is beneficial to introduce a prior term. \\  
MBIR methods allow the reformulation of the inverse problem as a minimization problem involving a fit-to-data function and a regularization (or prior) function.
In this paper, we consider the popular formulation given by:
\begin{equation}
 \min_{\bx\ge 0} ||A\bx - \by||_2^2 + \lambda  \mathcal{R}(\bx), 
 \label{eq:LS+Reg+constrain}
 \end{equation}
where the least squares term represents the data-fitting function, $\mathcal{R}$ is the prior term, and $\lambda$ is the regularization parameter.
Compressed sensing theory plays a crucial role in this framework, enabling the reconstruction of high-quality images from a limited amount of data through a suitable choice of the prior  $\mathcal{R(\bx)}$.
Specifically, forcing gradient sparsity through the prior has been proven to be extremely useful in fields such as medical imaging, where preserving boundaries of low-contrast objects is essential for accurate diagnosis and treatment.
Given a fixed $0 < p \le 1$, the gradient sparsity-inducing prior $\mathcal{R}$ is defined as: 
\begin{equation}\label{eq:definition_of_TpV}
    \mathcal{R}(\bx) =  || D \bx ||_{2, p}^p := || \>  | D \bx | \> ||_p^p,
\end{equation}
where  $| D \bx | \in \R^n$ represents the gradient-magnitude image of $\bx$, given by  $\left( | D \bx | \right)_i:= \sqrt{(D_h \bx)_i^2 + (D_v \bx)_i^2\ }$ with the operators $D_h$ and $D_v$, denoting the differences among two adjacent pixels in the horizontal and vertical directions, respectively.

It is widely acknowledged that the $\ell_0$ quasi-norm (given by Equation \eqref{eq:definition_of_TpV} when $p\to 0$), which quantifies the cardinality of the non-zero entries in its argument, represents the optimal sparsifying prior. 
Unfortunately, minimizations involving an $\ell_0$ prior pose significant computational challenges and are seldom addressed directly. 
Luckily, signal recovery remains feasible by substituting the $\ell_1$ norm (with $p=1$ in Equation \eqref{eq:definition_of_TpV}) for the $\ell_0$ quasi-norm under appropriate assumptions \cite{candes2006robust, donoho2006compressed}. In this case the regularization function is denoted as Total Variation (TV), which is very popular in CT applications. 
Many scholars also conducted a large amount of research on TV regularization, proposing numerous improved convex TV-based models, as we will detail in the next Section.  Nonetheless, the TV prior exhibits limitations such as its tendency to excessively smooth regions containing fine details or textures. 

Setting $0<p<1$ in Equation \eqref{eq:definition_of_TpV} defines $\mathcal{R}$ as the isotropic Total $p$-norm Variation (\TpV), which represents a better alternative to the TV prior for approximating the $\ell_0$ quasi-norm. 
The \TpV\ has already been successfully used in a number of imaging applications, including CT reconstruction, 
but its main drawback lies in its non-convex nature. Thus, ensuring the attainment of the global minimum by an MBIR may not be feasible, and MBIR approaches are prone to fall in one of the local minima, which could be significantly distant from the desired global solution. 
Consequently, non-convex minimization presents two practical disadvantages: the extensive computational effort required by the MBIR solver to converge and the unreliable dependence of the reconstructed image on the algorithm initialization.
Nevertheless, the identification of a \textit{good local} minimum can yield accurate reconstructions of sparse signals using significantly fewer measurements compared to the $p=1$ case \cite{chartrand2007exact}. \\
Given the urgent need for rapid and accurate images in clinical applications, Deep Learning (DL) has emerged as a powerful alternative to MBIR.
Many different approaches have been designed in the last years to address CT reconstruction with neural networks, and they have demonstrated the ability of data-driven methods to learn complex patterns from vast datasets, enabling the computation of high-quality images swiftly. 
However, three main drawbacks must be considered to make deep learning-based approaches usable in practice. 
Firstly, the black box structure of the trained networks makes it impossible to characterize the final outputs as numerical solutions of any mathematical inverse problem \cite{sidky2021docnn}.
Secondly, the training of the neural networks typically requires extensive training samples, which may not be available in many real-world medical scenarios. As we will discuss in Section \ref{sec:SOTA}, only a few works tackle this disadvantage in the literature.
Lastly, while highly accurate when used on images consistent with the training samples, neural networks may manifest untrustworthy effects when applied to slightly different data. This unstable behavior is the subject of study in some recent research \cite{antun2020instabilities, morotti2021green, evangelista2023ambiguity}.

\subsection*{Our proposal and goals} 

\begin{figure}
\centering
\includegraphics[width=0.95\textwidth]{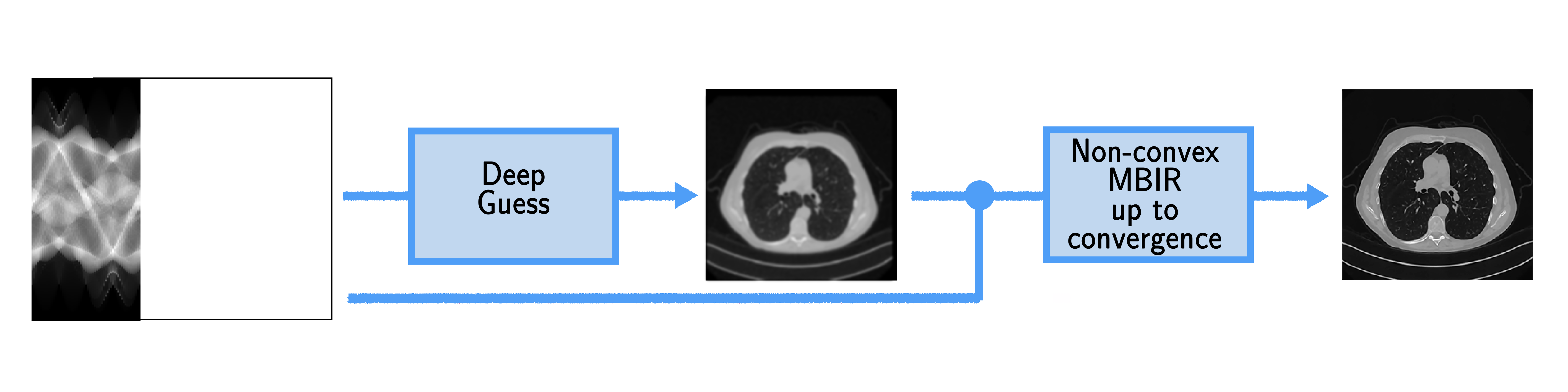}
\caption{Graphical representation of the proposed data-driven approach for solving the non-convex imaging problem.}
\label{fig:graph_abstr}
\end{figure}

This research proposes to the imaging community a fast and reliable tool for CT image reconstruction, taking into consideration the pros and cons of the aforementioned state-of-the-art approaches and directly addressing the scientific challenges for the implementation of software in clinical routines.\\
Specifically, we combine the advantages of both MBIR and DL-based best tools and design a unified solution of practical interest for medical imaging.
The core idea is to use a neural network to initialize the starting iterate in the algorithm addressing the \TpV-regularized problem. We refer to our strategy as the \textit{Deep Guess acceleration} approach to highlight that we exploit the effectiveness of deep learning to define the initial guess in an MBIR solver and make its numerical convergence faster.
The scheme of the Deep Guess acceleration framework is visually illustrated in Figure \ref{fig:graph_abstr}. The benefits of the proposed approach are manifold. 

\begin{itemize}
    \item Fast execution: by using a pre-trained neural network, the iterative solver can start from a position much closer to its solution, reducing the number of needed iterations.
    \item Reduced risk of bad local minima: with better initial guess, the likelihood of the solver falling into a bad image is significantly decreased.
    \item Mathematical explainability of the solution: running an MBIR to its convergence, we fully preserve convergence properties and the numerical characterization of the final solution, despite the use of a neural network.
\end{itemize}
At last, we highlight that we have also implemented neural networks using a ground truth-free supervised learning approach. Given the need for a network to be trained on a representative dataset, we propose incorporating into the Deep Guess step a network trained on images generated from available tomographic data. This offers a viable solution in scenarios where ground truth images are not available for standard network training.

The paper is organized as follows. In Section \ref{sec:SOTA} we present a state-of-art on both MBIR and DL-based algorithms for sparse-view CT image reconstruction; in Section \ref{sec:framework} we describe in detail the implementation of the proposed framework.
We report in Section \ref{sec:numerical_results} the results of several numerical experiments on both a real and a synthetic data set, demonstrating that our framework achieves good solutions, outperforming iterative non-convex algorithms and the  data-driven methods considered for comparison. Lastly, in Section \ref{sec:conclusion} we state some conclusions.

\section{State of the art \label{sec:SOTA}}

CT image reconstruction techniques have been evolving through two significant classes of tools: model-based iterative algorithms and deep learning-based tools. Both methodologies seek to enhance the precision and efficiency of computed tomography (CT) image reconstruction. MBIR  emphasizes iterative optimization processes and the incorporation of prior knowledge. In contrast, Deep Learning-based techniques leverage data-driven learning to achieve superior performance.
Here, we review some of the most impactful techniques proposed over the years. \\
The first class we focus on is the one composed of MBIR approaches, which iteratively solve optimization problems as the one considered in \eqref{eq:LS+Reg+constrain}. To produce high-quality images from sparse and noisy data, effectively reduce artifacts, and enhance edge preservation, MBIR schemes incorporate specific regularizers $\mathcal{R}(\bx)$ so that the image processing is performed according to precise mathematically-grounded image priors.
In medical imaging, priors that enforce gradient sparsity have become a well-established de facto choice because most parts of medical images are relatively smooth and uniform, with significant information concentrated at the edges. \\
Total Variation (TV) prior, minimizing the $\ell_1$ norm of the gradient-magnitude image, has been extensively utilized since its seminal paper by Rudin, Osher and Fatemi \cite{rudin1992nonlinear}, to effectively denoise images  (see \cite{vogel1996iterative, chan2005recent, rodriguez2013total}, to cite a few).
These formulations effectively differentiate between artifacts (small, random intensity changes) and significant image features (sharp edges), making TV prior interesting for tomographic applications and particularly beneficial in sparse data scenarios, where sparse-view artifacts and noise must be removed. This success has been proved in many works, such as \cite{persson2001total, sidky2008image, tian2011low, ritschl2011improved, caselles2015total, piccolomini2016fast, lv2020nonlocal, loli2021model}.\\
The TV models mentioned above are all convex and relatively easy to handle numerically. 
At the same time, it is well known that non-convex norms are more suitable for measuring sparsity than the corresponding convex norms, since they are much closer to the $\ell_0$ norm (that is exactly the measure of sparsity) than their convex counterparts \cite{chartrand2007exact, wang2018nonconvex}. 
For instance, the \TpV\ prior, which uses the $p$-norm of the gradient magnitude with $0<p<1$, allows for more flexible and adaptive regularization than TV.
It, in fact, tends to avoid staircasing effects better than TV, because it adjusts the degree of sparsity imposed on the gradients. For this reason it has been successfully used in medical imaging   \cite{tong2018edge, wang2021non} and for CT reconstruction  \cite{sidky2014cttpv, demirel2021p}. 
The wavelets have been used in tomographic imaging in \cite{borsdorf2008wavelet, beck2009fastwavelet, mehranian2013x, xu2016accelerated, purisha2017controlled}.\\
The second class of algorithms to tackle CT reconstruction is the one involving Deep Learning. 
DL-based methods offer faster reconstructions and can capture intricate patterns and structures, often outperforming traditional techniques in terms of speed and image quality. 
Nowadays, there is a plethora of papers exploring different ways in which neural networks can be used for medical image reconstruction. \\
One approach is to learn the reconstruction operator/solver, as the mapping between projection and image domain, from an appropriate corpus of training data \cite{zhu2018image, boink2019learned}.
A further approach implements two-step frameworks. In the first step, a reconstruction algorithm is applied offline to quickly compute the reconstruction from the projection data (typically, an analytic solver as the FBP), and then a trained neural network is applied as a post-processing operator to remove noise and artifacts \cite{jin2017deep,chen2017low, schnurr2019simulation, morotti2021green}. This approach is often referred to as Learnt Post Processing (LPP). 
 Both these approaches are typically very fast since they do not need the computation of the projector operator $A$ and any related matrix-vector multiplications (that are computationally heavy steps in MBIR), but the mathematical characterization of the solution is questionable \cite{sidky2021docnn}. 
Alternatively, there are hybrid approaches where DL and MBIR are merged together. 
Examples are the unrolling algorithms, which implement iterative solvers (as classical optimization algorithms) unrolling each iteration in a layer of a deep neural network \cite{monga2021algorithm, xiang2021fista}. 
Here, the optimal model parameters are learned during the training phase so that adaptive parameters allow high-quality reconstructions, and the convergence of the MBIR solver is accelerated too. 
A popular option is represented by the class of Plug-and-Play (PnP) hybrid methods, where a learned regularizer is plugged inside a classical variational scheme \cite{venkatakrishnan2013plug, cascarano2022plug}. 
Specifically, the network plays the role of a denoiser in one inner step of each iteration of the MBIR, which is executed in a classic way up to its convergence.
A further hybrid approach is represented by Deeply Learned Regularizers, which employ a pre-trained neural network to learn the  $\mathcal{R}$ prior operator. An example is represented by the NETT algorithm \cite{li2020nett}.
An additional class of hybrid methods is represented by multi-step schemes. To set some examples, in \cite{evangelista2023rising},  the MBIR method is incorporated into a ground truth-free Learned Post Processing scheme. Here, a coarse image is initially reconstructed from projections using an early interrupted MBIR solver and then processed by a network trained on accurate MBIR reconstructions.  In a second example \cite{morotti2024space}, a neural network adjusts an adaptive space-variant total variation (TV)-based regularizer, and the resulting model-based formulation is subsequently solved efficiently.\\
The aforementioned papers use supervised learning, heavily relying on the availability of a set of consistent and clean ground truth images, which is not always guaranteed in real medical scenarios. 
To the best of our knowledge, few researches have tackled it so far.
The authors of \cite{baguer2020computed, shu2022sparse} propose to use the Deep Image Prior approach, which is a scheme where the architecture of an untrained neural network itself serves as a regularizer for imaging tasks. The network is initialized with random weights and iteratively fitted to the observed datum $\by$  so that the network naturally captures the image's features and promotes accurate reconstructions.
Alternatively, in \cite{superultra}, an unsupervised hybrid reconstruction is coupled with a supervised denoiser, which does not need medical ground truth images, whereas in \cite{evangelista2023rising}, the authors generate the target images for the network training from the available set of projection data, by running MBIR accurate reconstructions offline.

\section{The proposed Deep Guess acceleration framework \label{sec:framework}}

As anticipated in the Introduction, we consider in this paper a model-based approach formulated as the minimization problem \eqref{eq:LS+Reg+constrain} with regularization prior given by \eqref{eq:definition_of_TpV}. We focus on the non-convex case with $0<p<1$, leading to the so called \TpV\ regularization.\\
The proposed  framework comprises two primary steps, as illustrated in Figure \ref{fig:graph_abstr}.
\begin{enumerate}
\item 
The initial step involves calculating a robust approximation of the reconstructed image using a neural network within a supervised  post processing scheme. The output of this step, referred to as Deep Guess (DG), is a reconstruction scheme already proposed in the literature (see Section \ref{sec:SOTA}). Here, however, it serves as the input for the subsequent phase.
\item 
The second step refines the quality of image computed by the DG by applying to it some iterations of an iterative optimization algorithm solving the non-convex  minimization   problem \eqref{eq:LS+Reg+constrain}.
\end{enumerate}

We now describe in detail the implementation of the two steps, starting from the second one.

\subsection{A non-convex Model-Based Iterative Reconstruction method utilizing Chambolle-Pock iterations \label{subsec:nonconvex MBIR}}

The  minimization \eqref{eq:LS+Reg+constrain} with non-convex \TpV\ prior is solved by 
the iterative reweighting $\ell_1$-norm (IR$\ell_1$) strategy introduced in \cite{candes2008enhancing,daubechies2010iteratively}. It is an iterative process, where each iterate is obtained by solving a convex optimization problem, i.e.:
\begin{equation}\label{eq:inverse_problem_weighted}
     \bx^{(k+1)} = \arg\min_{\bx\ge0} \| A\bx - \by \|_2^2  + \lambda || \ \bw(\bx^{(k)}) \odot | D\bx| \ ||_1,
\end{equation}
where the weights $\bw(\bx^{(k)}) \in \R^n$ are defined as:
$$
\bw(\bx^{(k)})_i = \Bigg( \frac{ \sqrt{\eta^2 + | D \bx^{(k)}|_i^2}}{ \eta } \Bigg)^{p-1},
$$
with a predefined smoothing parameter
$\eta>0$.\\
Each problem \eqref{eq:inverse_problem_weighted} is approximately solved using a limited number of iterations of the Chambolle-Pock (CP) algorithm, a primal-dual method originally introduced for the minimization of convex functional appearing in image processing applications in \cite{chambolle2011first}. The CP algorithm has been effectively employed for total variation (TV)-based image reconstruction from undersampled tomographic data in \cite{sidky2012convex, loli2021model}. In \cite{sidky2014cttpv}, the authors suggest applying the CP algorithm in a manner consistent with the approach used in this study.\\
In its original formulation, the CP algorithm was proposed to minimize an objective function of the form:
\begin{equation}\label{eq:CP_general_formulation}
    \min_{\bx \in \R^n} F(M\bx) + G(\bx),
\end{equation}
where both $F$ and $G$ are real-valued, proper, convex, lower semi-continuous functions, and $M$ is a linear operator from $\R ^n$ to $\R^s$. 
Since there are no constraints on the smoothness of either $F$ and $G$, the CP method can be applied to our problem by setting:
\begin{equation*}\label{eq:CP_funtions}
    \begin{cases}
        G(\bx) = \iota_{\mathcal{X}}(\bx), \\
        F(M\bx) = \frac{1}{2} || A\bx - \by ||_2^2 + \lambda || \ \bw \odot | D\bx| \ ||_1,
    \end{cases}    
\end{equation*}
where $\iota_{\mathcal{X}}(\bx)$ is the indicator function of the feasible set $\mathcal{X} := \{\bx \in \R^n; \bx_i \geq \boldsymbol{0} \ \forall i=1, \dots n\}$, and the linear operator $M \in \R^{s \times n}$ is given by the row-wise concatenation of $A$ and $D$, such that $M = \left[ A; D \right]$. 
To preserve a smooth and readable flow, we omit some detailed calculations that can be found in \cite{morotti2024space}. We only remark that we need to define $\bz^{(k)} = (p^{(k)}, q^{(k)})$ with two dual variables $p \in \R^m$ and $q \in \R^{2n}$,  to explicitly derive the $\prox_{\sigma F^*}$ term, so that the CP algorithm is defined by the following update rules:
\begin{align}
    \begin{cases}
        p^{(k+1)} = \frac{p^{(k)} + \sigma \left(A \bar{\bx}^{(k)} - \by^\delta \right)}{1 + 3\sigma}, \\
        q^{(k+1)} = \frac{\lambda \bw(\bar{\bx}^{(k)}) \odot \left(q^{(k)} + \sigma | D \bar{\bx}^{(k)} | \right)}{\max \left( \lambda \bw(\bar{\bx}^{(k)}),\ q^{(k)} + \sigma | D \bar{\bx}^{(k)} | \right)},\\
        \bx^{(k+1)} = \mathcal{P}_+\left(\bar{\bx}^{(k)} - \tau M^T\begin{bmatrix}
            p^{(k+1)} \\ q^{(k+1)}
        \end{bmatrix} \right), \\
        \bar{\bx}^{(k+1)} = \bx^{(k+1)} + \alpha \left(\bx^{(k+1)} - \bx^{(k)} \right).
    \end{cases}
\end{align}

In our numerical experiments, the CP iterations terminate when the distance between two consecutive iterates is less than a fixed tolerance $\tau$, i.e.:
$$ ||{\bx}^{(k+1)} - {\bx}^{(k)}||_2 <  \tau,$$ 
or when a maximum number of iterations is reached. The considered stopping criterion is frequently employed in image processing, where the images become indistinguishable if the tolerance is sufficiently small.

\subsection{The Deep Guess step \label{subsec:deep_guess}}

The Deep Guess  step  essentially uses a convolutional neural network as a post-processing step to refine a coarse reconstruction obtained by a fast algorithm, following a supervised learning approach.

\begin{figure}
\centering
\includegraphics[width=0.9\textwidth]{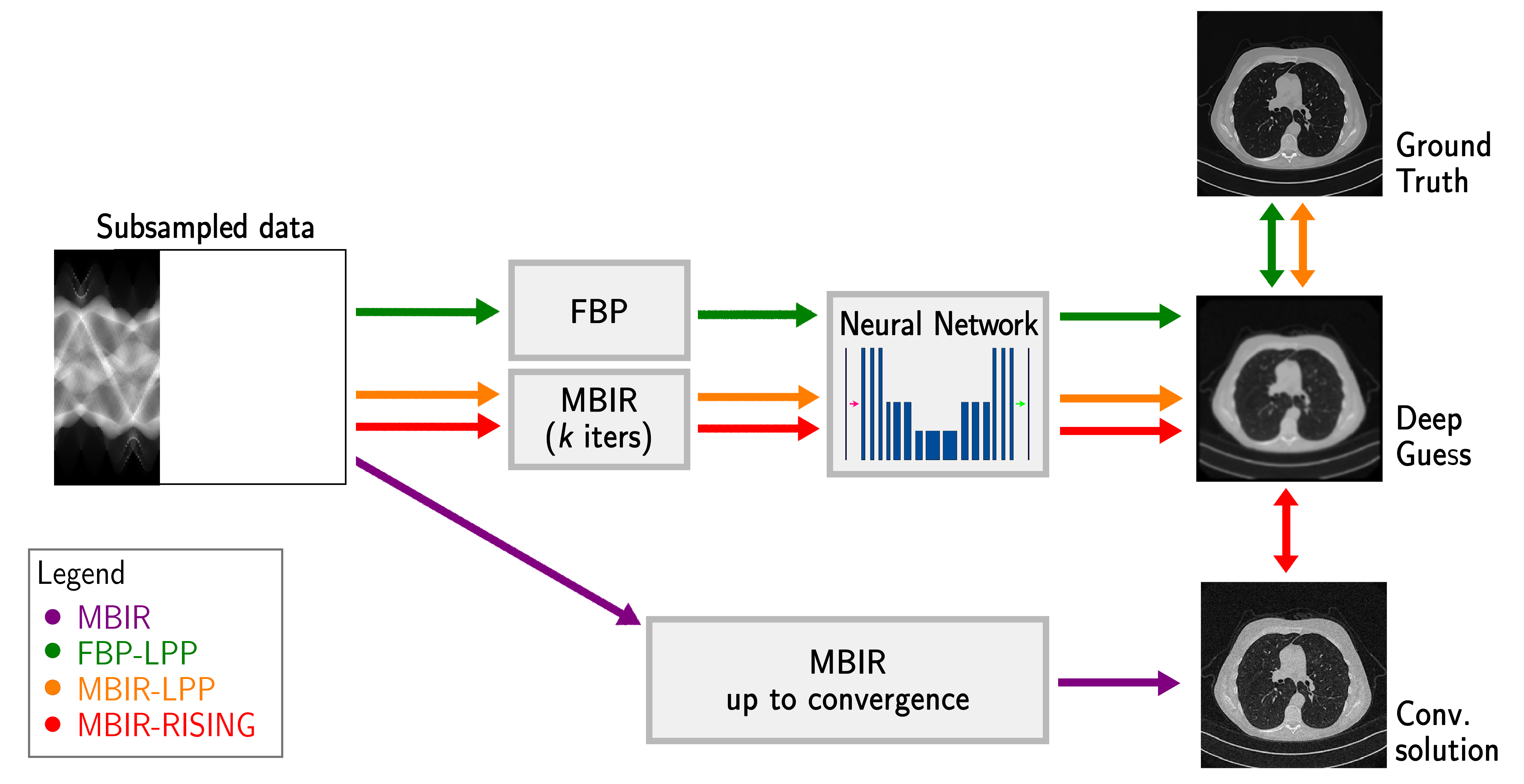}
\caption{The Deep Guess block. It can be computed by following the red, orange, or green path, all exploiting a neural network. The input to the network is a coarse reconstructed image, obtained either through FBP or after $k$ iterations of an MBIR solver. The target images are either the ground truth or those obtained when the MBIR method has converged.}
\label{fig:graph_abstr_DG}
\end{figure}

Figure \ref{fig:graph_abstr_DG} illustrates three distinct implementations of the Deep Guess step, represented by the red, orange, and green paths. Although all implementations utilize a neural network, they vary in the input and/or target images used. \\
The coarse input image of the neural network is derived from the sinogram either through a Filtered Back Projection (FBP) (green arrow) or through a limited number $K$ of iterations of an MBIR method (orange and red arrows).
As target images for the network training, the ground truth ones can be used if available.
In the absence of ground truth images, which is common in real medical imaging, we implement the RISING scheme proposed in \cite{evangelista2023rising}. In this case (red double-headed arrow), the targets are the images generated by the same MBIR method used  to compute the input images but now run until convergence (purple path). \\
Regarding the MBIR method used in the DG step, we emphasize that it is independent of the one used in the final step of Figure \ref{fig:graph_abstr}. 
In this paper we have tested the two following  MBIR methods to compute the Deep Guess.
\begin{itemize}
\item A Scaled Gradient Projection (SGP) method \cite{bonettini2008scaled, porta2015new,piccolomini2018reconstruction} applied to the convex minimization problem: 
 \begin{equation}
  \min_{\bx\ge0} \|A\bx - \by \|_2^2 + \mu TV_{\beta}(\bx)
  \label{eq:TV}
\end{equation}
where $TV_{\beta}$ is the differentiable TV prior defined as:
\begin{equation}\label{eq:definition_of_TV_beta}
     TV_{\beta}(\bx) :=  \sum_{i=1}^n \sqrt{\left( D_h \bx \right)_i^2 + \left( D_v \bx \right)_i^2+\beta^2},
 \end{equation}  
 for some small scalar $\beta > 0$ and $\mu$ is a positive regularization parameter.
 \item The  IR$\ell_1$ method utilizing CP iterations described in Section \ref{subsec:nonconvex MBIR} applied to the non-convex \TpV\ minimization:
 \begin{equation}
 \min_{\bx\ge 0} ||A\bx - \by||_2^2 + \mu || D \bx ||_{2, p}^p, 
 \label{eq:CP+TpV}
 \end{equation}
where the regularization parameter is denoted by $\mu$ to differentiate it from the parameter $\lambda$ used in the final step of the framework.
 \end{itemize}

As a neural network for the image-to-image step, we have used an architecture inspired by a classical UNet \cite{ronneberger_2015_unet}. Specifically, we utilize the residual UNet (ResUNet) architecture, which has been previously employed for image reconstruction from sparse-view CT in \cite{morotti2021green, evangelista2023rising}.
The network is trained using a loss function that represents the mean squared error  between the network predictions and the corresponding target values. The Adam optimizer is utilized with a learning rate of $0.001$, and the training process is conducted over a total of 50 epochs.
It is noteworthy that the forward pass of a pre-trained neural network is characterized by very short computational times. Consequently, the entire Deep Guess step can be executed rapidly.
At last, we remark again the described post-processing approaches are not new in the literature, but so far, they have been used as solvers for the CT inverse problem. With their usage in the Deep Guess Acceleration scheme, we leverage the state-of-the-art literature strongly relying on deep learning as an efficient tool inside a more interpretable mathematical framework, providing a variational solution as the final reconstruction.

\section{Numerical experiments and discussion}\label{sec:numerical_results}

In this Section we present  the numerical experiments performed on  test problems with both real and synthetic medical images. 
The code to replicate the experiments can be found at: \url{https://github.com/devangelista2/DeepGuess}.\\
Here, we first describe the datasets and test problems utilized in our experiments. Subsequently, we analyze the results obtained on real medical images using the proposed approach, comparing them to other methods. Finally, we conduct a more detailed examination of the possible implementations of the Deep Guess step and its stability properties. \\
In this Section, we will use the following notation to indicate the different Deep Guess implementations. 
According to the notations used in Figure \ref{fig:graph_abstr_DG}, each version is designated by a two-part acronym. The first part specifies the algorithm used to compute the network's input, and the second part, separated by a hyphen, indicates the type of target used during network training: `LPP' if the target is a ground truth, and `RISING' if the target is an image  generated by running the MBIR model to convergence.\\
We also clarify that the number of iterations reported for the CP method represents the total number of iterations carried out during the IR$\ell_1$ scheme, calculated as the sum of CP iterations across all iterations of the IR$\ell_1$ procedure.

The metrics used for comparison are the Relative Error (RE) between the computed reconstruction and the ground truth and the Structural Similarity Index (SSIM) \cite{wang2003multiscale}, whose values are scaled in the range $[0,100]$.\\

\subsection{The datasets \label{subsec:datasets}}

We considered both a dataset of real CT images and synthetic images. 

\begin{figure}[H]
\begin{minipage}[outer sep=0]{\textwidth}
\centering
\begin{minipage}[m]{0.23\textwidth}
    \begin{tikzpicture}
        \node [anchor=south west, inner sep=0] (image) at (0,0) {\includegraphics[width =\textwidth]{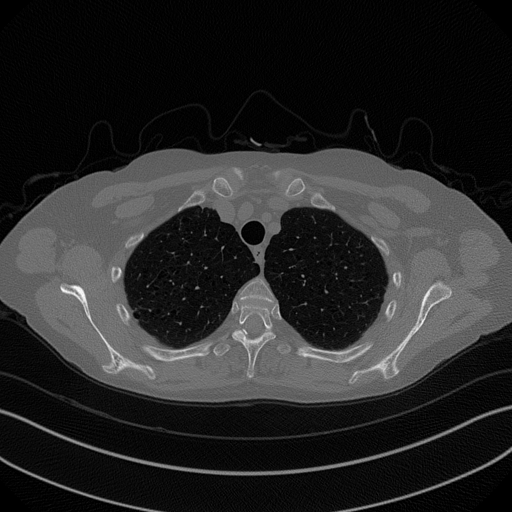}};
        \begin{scope}[x={(image.south east)}, y={(image.north west)}]
            \draw[red, thick] (0.42, 0.22) rectangle (0.85, 0.44);
        \end{scope}
        \begin{scope}[x={(image.south east)}, y={(image.north west)}]
            \draw[red, thick] (0.14, 0.48) rectangle (0.57, 0.69);
        \end{scope}
    \end{tikzpicture}
    \end{minipage} 
\begin{minipage}[m]{0.22\textwidth}
    \includegraphics[trim= 20mm 63mm 55mm 40mm,clip,width=\textwidth]{imm/mayo107_GT.png}\\[1ex]
    \includegraphics[trim=55mm 30mm 20mm 73mm, clip, width=\textwidth]{imm/mayo107_GT.png}
    \end{minipage}  
\quad
\begin{minipage}[m]{0.23\textwidth}
\begin{tikzpicture}
        \node [anchor=south west, inner sep=0] (image) at (0,0) {\includegraphics[width =\textwidth]{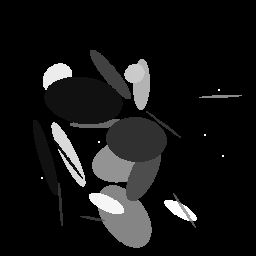}};
        \begin{scope}[x={(image.south east)}, y={(image.north west)}]
            \draw[red, thick] (0.34, 0.09) rectangle (0.85, 0.30);
        \end{scope}
        \begin{scope}[x={(image.south east)}, y={(image.north west)}]
            \draw[red, thick] (0.07, 0.38) rectangle (0.55, 0.60);
        \end{scope}
    \end{tikzpicture} 
    \end{minipage}
\begin{minipage}[m]{0.22\textwidth}
    \includegraphics[trim= 5mm 25mm 30mm 25mm,clip,width=\textwidth]{imm/coule14_GT.png}\\[1ex]
    \includegraphics[trim=22mm 5mm 13mm 45mm, clip, width=\textwidth]{imm/coule14_GT.png}
    \end{minipage}
\end{minipage}
    \caption{\label{fig:gt} Ground truth images from the test set of Mayo Clinic (left) and COULE (right), with two zooms-in on the regions depicted with red rectangles on the corresponding whole images. }
\end{figure}

For each ground truth image $\bx^{GT}$ in a dataset, we sampled a noise realization $\be \sim \mathcal{N}(\boldsymbol{0}, I)$, and we computed:
\begin{equation}\label{eq:pbtest}
\by=A \bx^{GT} + \nu \frac{||A \bx^{GT}||}{||\be||} \be,
\end{equation}
where $\nu \geq 0 $ is the noise level.
The matrix $A$ has been computed with Astra Toolbox \cite{van2015astra,van2016fast} functions, simulating a sparse fan beam geometry with $n_{angle}$ projections uniformly distributed within the angular range $(0,n_{\text{range}})$ degrees, with a detector resolution of $2 \sqrt{n}$ pixels. In the following, we indicate this setup as the tomographic geometry $\mathcal{G}_{n_{range},n_{angle}}$.\\
The dataset of real CT images considered is the Mayo Clinic dataset \cite{mccollough2016tu}, constituted of $512 \times 512$ pixel real images of human abdomen. We used 3306 images for training and 357 for testing.  We generated projections as in \eqref{eq:pbtest} on a flat detector composed of 1024 pixels, using different geometries. In all the experiments, we added noise with $\nu=0.001$. 
An example image is presented in Figure \ref{fig:gt}, accompanied by two zoomed-in views of regions of interest.   \\
The considered synthetic data set is the Contrasted Overlapping Uniform Lines and Ellipses (COULE) dataset, downloadable from 
\url{www.kaggle.com/loiboresearchgroup/coule-dataset}. It contains 430  synthetic sparse-gradient gray-scale images of size $256 \times 256$ with many overlying objects, varying in size and contrast with respect to the background. Synthetic images are ideal representations of objects, free from any noise or artifacts. An example of a COULE image is displayed in Figure \ref{fig:gt}, together with two zoomed crops of interest.  We generated projections as in \eqref{eq:pbtest} on a flat detector composed of 512 pixels, using the geometry $\mathcal{G}_{180,180}$. Noise has been  added with  $\nu=0.01$. We divided the dataset into a training set comprising 400 images and a test set containing 30 images.  \\
Both with the Mayo Clinic and the COULE datasets, we used the available ground truth images for training the neural network and getting the LPP Deep Guess images (corresponding to the orange and green final arrows, in Figure \ref{fig:graph_abstr_DG}). We also generated the Mayo Clinic and COULE datasets of MBIR convergent solutions, to train the DG network with the RISING principle (corresponding to the final red arrow in Figure \ref{fig:graph_abstr_DG}).  In this case, the regularization parameters have been set by trial and error on the training samples. 

\subsection{Preliminary  results with different MBIR approaches \label{subsec:variational}}

\begin{figure}
TV (SGP)  \hspace{35mm}  \TpV\ (CP)     \hspace{36mm}  W (FISTA)  \\    \includegraphics[width=0.3\textwidth]{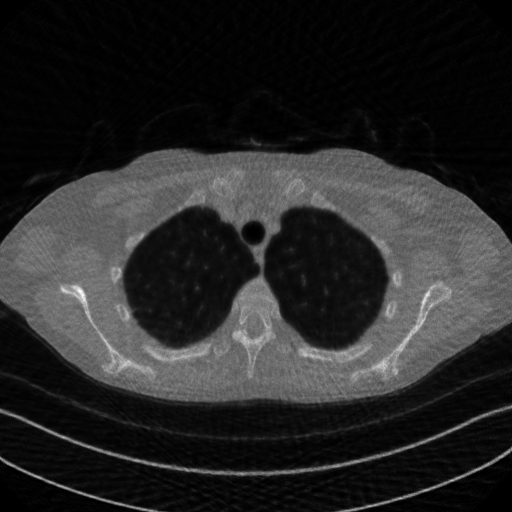}
    \includegraphics[width=0.3\textwidth]{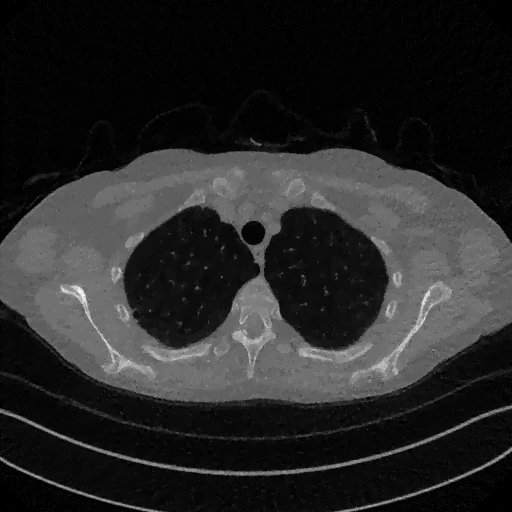}
    \includegraphics[width=0.3\textwidth]{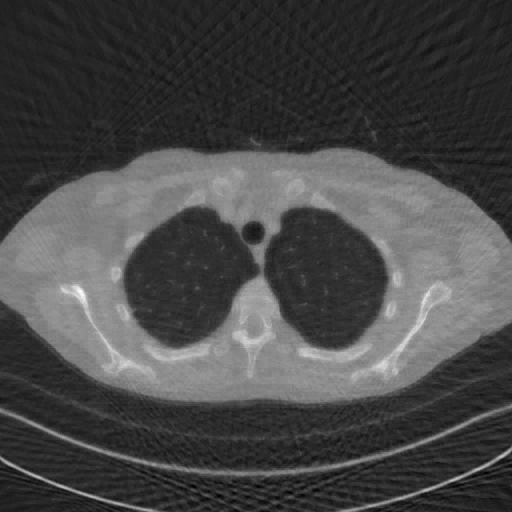}
    \\
    \caption{Results on the Mayo  test image shown in Figure \ref{fig:gt}, computed using the TV prior (first column), the \TpV\ prior (second column), and the W regularizer (third column), in case of $\mathcal{G}_{180,60}$ sparse protocol. }
    \label{fig:mayo_variazionali}
\end{figure}

In our preliminary experiments, we addressed the CT image reconstruction problem using various MBIR approaches and different geometries to validate the selection of the \TpV\ model as the final MBIR solver employed in all subsequent tests.\\
The first considered MBIR approach is  the Scaled Gradient Projection (SGP) algorithm applied to a $TV_{\beta}$-regularized model (with $\beta=10^{-3}$), as previously discussed in Section \ref{subsec:deep_guess}. It is named as TV(SGP) in the following.\\
The second alternative method is the Fast Iterative Shrink-Thresholding Algorithm (FISTA) applied to compute the optimal sparse wavelet (W) coefficients by solving the following minimization problem \cite{beck2009fast,beck2009fastwavelet}:
\begin{equation}
    \min_{\boldsymbol{c} \in \R^n} ||\bar{A}\boldsymbol{c} - \by||_2^2 + \lambda  ||\boldsymbol{c}||_1,
\end{equation}
where the matrix $\bar{A}=A \mathcal{W}^T$ is given by the product of the projection matrix $A$ and of the matrix $\mathcal{W}^T \in \R^{n \times n}$, representing the synthesis operator of a wavelet transform, and $\boldsymbol{c}$ represents the wavelet coefficients of the solution $\bx = \mathcal{W}^T\boldsymbol{c}$. In the following, we will denote as W(FISTA) this approach. \\
The results obtained on the test image previously displayed are shown in Table \ref{tab:imm107mayo_G360360_variaz}, where we also report for comparison the metrics of the traditional FBP method.
From the table, it is immediately apparent that FBP performs very poorly, particularly for the sparsest geometry, denoting the necessity for MBIR approaches in these scenarios.
In the settings of $\mathcal{G}_{360,360}$ and $\mathcal{G}_{180,180}$, the three compared MBIR methods exhibit very similar performance. However, in cases with very few views, the \TpV\ model significantly outperforms the others, exceeding the TV model by approximately $10\%$ and the $\ell_1$ wavelet by about $30\%$ in terms of both RE and SSIM.\\
Moreover, it is notable that all the MBIR approaches are quite slow, necessitating hundreds of iterations to achieve the final reconstructions.
Therefore, to make these methods applicable in real-world cases, it becomes of primary importance to reduce the number of iterations and, consequently, the computational time.\\
Figure \ref{fig:mayo_variazionali} displays the reconstructions of the chosen image obtained in the quite sparse case $\mathcal{G}_{180,60}$ and computed, from left to right, by the TV, \TpV\ and W models. The \TpV\ image is the most accurate and clean, whereas artifacts are clearly visible in the TV and wavelet ones, according to the metrics in Table \ref{fig:mayo_variazionali}.

\begin{table}[]
    \centering
     \setlength\tabcolsep{0pt}
    \begin{tabular*}{\linewidth}{@{\extracolsep{\fill}} ll  cc ccc }
        \toprule
  &     &  SSIM & RE     & iters   \\

\midrule    
\multirowcell{3}{$\mathcal{G}_{360,360}$ }
& \quad FBP             & 84.47 & 0.0744 &  -     \\ 
& \quad TV (SGP, $\lambda = 0.001$)        & 93.54 & 0.0469 & 150   \\
& \quad \TpV\ (CP, $\lambda = 0.0001$)        & 91.72 & 0.0539 & 500   \\
& \quad W (FISTA, $\lambda = 0.001$)       & 91.52 & 0.0556 &  500  \\
\midrule
\multirowcell{3}{$\mathcal{G}_{180,180}$ }
& \quad FBP             & 58.84 & 0.1466 &  -     \\
& \quad TV (SGP, $\lambda = 0.001$)        & 87.99 & 0.0662 & 300   \\
& \quad \TpV\ (CP, $\lambda = 0.0001$)        & 86.60 & 0.0699 & 500  \\
& \quad W (FISTA, $\lambda = 0.001$)       & 86.43 & 0.0726 &  500 \\
\midrule
\multirowcell{3}{$\mathcal{G}_{180,120}$ }
& \quad FBP             & 47.99 & 0.1855 &  -     \\
& \quad TV (SGP, $\lambda = 0.001$)        & 83.43 & 0.0835 & 300   \\
& \quad \TpV\ (CP, $\lambda = 0.0001$)        & 83.95 & 0.0803 & 500  \\
& \quad W (FISTA, $\lambda = 0.001$)       & 79.05 & 0.0981 &  500  \\
\midrule
\multirowcell{3}{$\mathcal{G}_{180,60}$ }
& \quad FBP             & 30.56  & 0.2955  &  -     \\
& \quad TV (SGP, $\lambda = 0.001$)        & 74.81 & 0.1190 & 300   \\
& \quad \TpV\ (CP, $\lambda = 0.0001$)        & 78.08 & 0.1039 & 500  \\
& \quad W (FISTA, $\lambda = 0.001$)       & 63.25 & 0.1517 & 500 \\
\midrule
\multirowcell{3}{$\mathcal{G}_{180,30}$ }
& \quad FBP             & 18.12 & 0.4698 &  -     \\
& \quad TV (SGP, $\lambda = 0.01$)        & 66.65 & 0.1646 & 300 \\
& \quad \TpV\ (CP, $\lambda = 0.0001$)        & 71.11 & 0.1464 & 500   \\
& \quad W (FISTA, $\lambda = 0.01$)       & 51.70 & 0.2093 &  500 \\
\bottomrule
    \end{tabular*}
    \caption{Values of the quantitative metrics, computed on the Mayo Clinic test image reconstructions, for different MBIR methods.  }
    \label{tab:imm107mayo_G360360_variaz}
\end{table}

\newpage
\subsection{Results with the proposed Deep Guess Acceleration framework\label{subsec:mayo}}

\begin{figure}
CP  \hspace{45 mm} CP    \hspace{45mm} CP   \\ 
from zeros\hspace{35 mm} from DG by FBP-LPP     \hspace{17mm} from DG by TV-RISING   \\
    \includegraphics[width=0.3\textwidth]{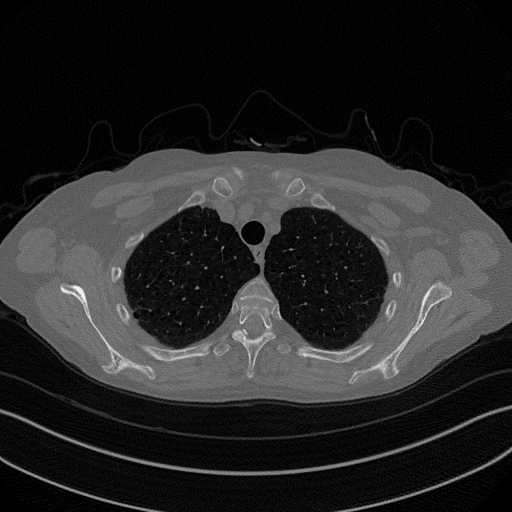}
    \includegraphics[width=0.3\textwidth]{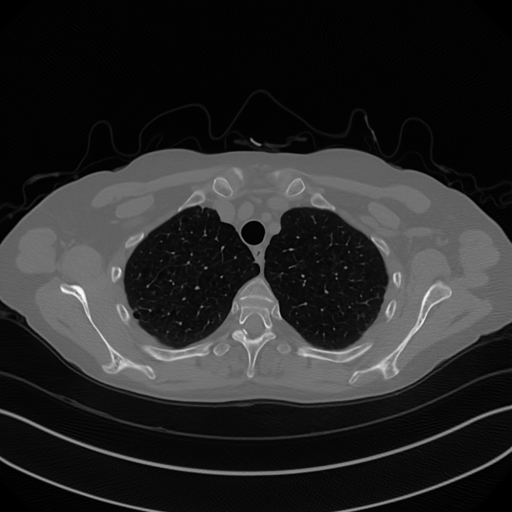}
    \includegraphics[width=0.3\textwidth]{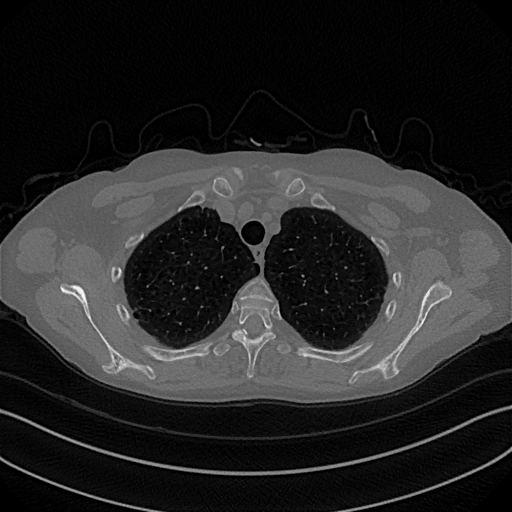}
    \\
    \includegraphics[trim=20mm 63mm 55mm 40mm, clip, width=0.3\textwidth]{imm/mayo107_360_zeros_IS.png}
    \includegraphics[trim=20mm 63mm 55mm 40mm, clip, width=0.3\textwidth]{imm/mayo107_360_FbpLpp_IS.png}
    \includegraphics[trim=20mm 63mm 55mm 40mm, clip, width=0.3\textwidth]{imm/mayo107_360_SgpRising_IS.png}\\
    \includegraphics[trim=55mm 30mm 20mm 73mm, clip, width=0.3\textwidth]{imm/mayo107_360_zeros_IS.png}
    \includegraphics[trim=55mm 30mm 20mm 73mm, clip, width=0.3\textwidth]{imm/mayo107_360_FbpLpp_IS.png}
    \includegraphics[trim=55mm 30mm 20mm 73mm, clip, width=0.3\textwidth]{imm/mayo107_360_SgpRising_IS.png}
    \caption{Results on the Mayo Clinic test image shown in Figure \ref{fig:gt}, computed using the CP algorithm (first column) and the proposed FBP-LPP and TV-RISING Deep Guess approaches (second and third columns), within the $\mathcal{G}_{360,360}$ tomographic protocol and setting $p=0.5$ for the final CP execution.}
    \label{fig:mayo360}
\end{figure}

\begin{figure}
CP  \hspace{45 mm} CP    \hspace{45mm} CP   \\ 
from zeros\hspace{35 mm} from DG by FBP-LPP     \hspace{17mm} from DG by TV-RISING   \\
    \includegraphics[width=0.3\textwidth]{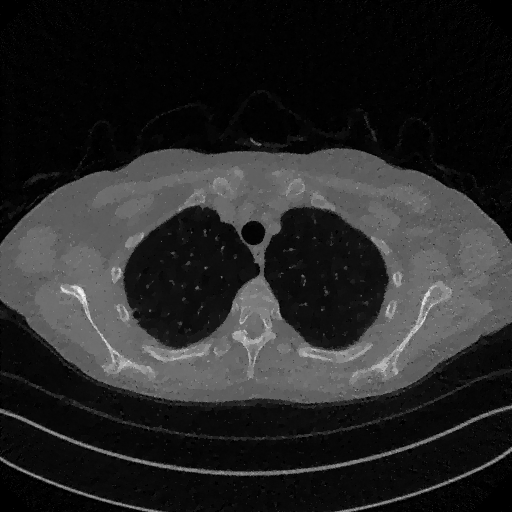}
    \includegraphics[width=0.3\textwidth]{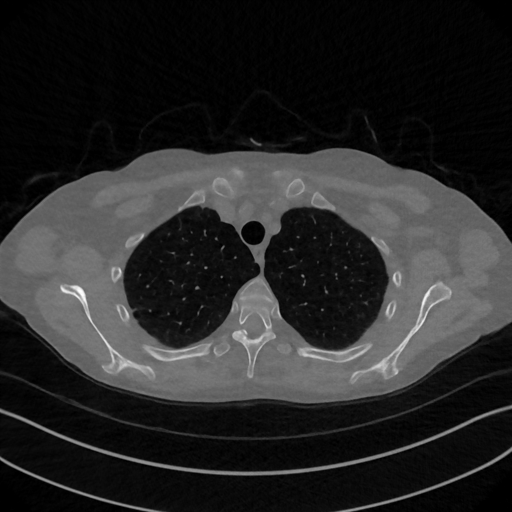} 
    \includegraphics[width=0.3\textwidth]{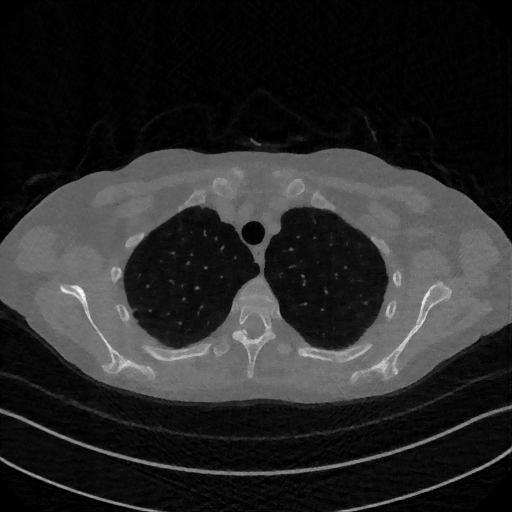} \\
\includegraphics[trim=20mm 63mm 55mm 40mm, clip, width=0.3\textwidth]{imm/mayo18060_p02_zeros_IS_E.png}
    \includegraphics[trim=20mm 63mm 55mm 40mm, clip, width=0.3\textwidth]{imm/mayo18060_p02_FbpLpp_IS.png} 
    \includegraphics[trim=20mm 63mm 55mm 40mm, clip, width=0.3\textwidth]{imm/mayo18060_p02_SgpRising_IS_36iter.png} \\ 
\includegraphics[trim=55mm 30mm 20mm 73mm, clip, width=0.3\textwidth]{imm/mayo18060_p02_zeros_IS_E.png}
    \includegraphics[trim=55mm 30mm 20mm 73mm, clip, width=0.3\textwidth]{imm/mayo18060_p02_FbpLpp_IS.png} 
    \includegraphics[trim=55mm 30mm 20mm 73mm, clip, width=0.3\textwidth]{imm/mayo18060_p02_SgpRising_IS_36iter.png} \\ \vspace{0mm}
\hspace{50mm} DG by FBP-LPP  \hspace{25mm} DG by TV-RISING  \\ \vspace{0mm}
    \hspace{50mm}
    \includegraphics[trim=20mm 63mm 55mm 40mm, clip, width=0.3\textwidth]{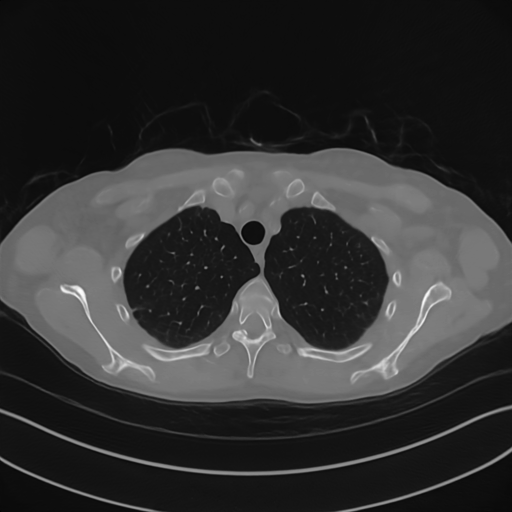} 
    \includegraphics[trim=20mm 63mm 55mm 40mm, clip, width=0.3\textwidth]{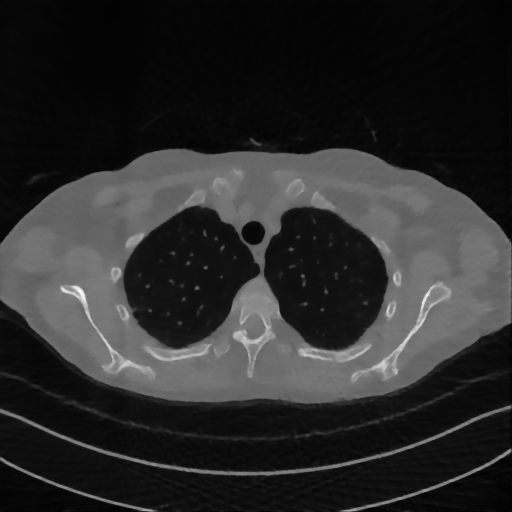} \\ \vspace{0mm} 
\hspace{50mm} 
    \includegraphics[trim=55mm 30mm 20mm 73mm, clip, width=0.3\textwidth]{imm/mayo18060_p02_FbpLpp_DIG.png} 
    \includegraphics[trim=55mm 30mm 20mm 73mm, clip, width=0.3\textwidth]{imm/mayo18060_p02_SgpRising_DIG.png} \\
\caption{Results on the Mayo Clinic test image, computed using the CP algorithm with different  initial guess: null (first column) and the proposed FBP-LPP and TV-RISING Deep Guess approaches (second and third columns), within the sparse $\mathcal{G}_{180,60}$ tomographic protocol and setting $p=0.2$ for the final CP execution. The last row displays the zooms of the images generated by the corresponding Deep Guess.}
    \label{fig:mayo60}
\end{figure}

In this Section, we present the results obtained by the proposed Deep Guess accelerated framework on the Mayo Clinic test image.\\
When the input of the neural network in the DG step is the image computed by $K$ MBIR iterations (red and orange arrows in Figure \ref{fig:graph_abstr_DG}, we set $K=10$ for $\mathcal{G}_{360,360}$ and $K=15$ for $\mathcal{G}_{180,60}$. 
Concerning the regularization parameter $\mu$ in \eqref{eq:TV} and \eqref{eq:CP+TpV}, it has been set by trial and error to make the CP perform at its best when $\bx^{(0)}$ is the null image, and it never changes when a different initial guess is used. \\
 The comparison includes results obtained using different initial guesses $\bx^{(0)}$ for the  Chambolle-Pock method, specifically a null vector, the analytical FBP solution, $K$ iterations of SGP applied to the SGP(TV), and two implementations of the Deep Guess.
Additionally, we compare our DG Accelerated framework with state-of-the-art data-driven approaches. Among the various methods available, we consider the FBP-LPP and TV-RISING approaches, here evaluated also as final reconstructors. We also include the NETT algorithm \cite{li2020nett}, popular for its mathematically explainable model.\\  
Table \ref{tab:imm107mayo_G360360} presents the error metrics for the test image, on the geometry $\mathcal{G}_{360,360}$. We first observe that for all the choices of the initial guess (except for the worst case when $\bx^{(0)}$ is computed by the TV model) the CP method produces similar results in terms of accuracy. It also requires significantly many iterations to converge when it starts from the zero-image or from the FBP or the TV-based solutions, whereas the starting guesses provided by the two considered data-driven strategies make the CP converge extremely fast.
Notably, in 29 and 66 iterations, the CP solver improves the FBP-LPP and the TV-RISING results, respectively, meaning that our \TpV-based MBIR overcomes the two state-of-the-art hybrid approaches. The NETT algorithm performs quite well, but its metrics are inferior compared to those of the CP algorithm.
Figure \ref{fig:mayo360} provides a visual representation of the results, confirming the high quality of all the displayed reconstructions which preserve all the thin details present in the ground truth, coherently with the achieved good metrics.\\
Focusing on the experiments performed with sparser tomographic acquisitions, we report the metrics relative to the $\mathcal{G}_{180,60}$ geometry in Table \ref{tab:imm107mayo_G18060} obtained by TpV with both $p=0.5$ and $p=0.2$.
At first glance, the value of $p$ does not significantly affect the results. In all instances, the proposed framework greatly accelerates the CP method compared to initial guesses of zeros, FBP, or TV, reducing computational time by 200-400 times. Additionally, it enhances the accuracy of our comparing methods FBP-LPP, TV-RISING, and NETT. Notably, the FBP-LPP implementation of Deep Guess performs the best.
As shown in Figure \ref{fig:mayo60} depicting the results by $p=0.2$ for this very sparse geometry,
 the CP reconstruction obtained with null initial guess starting still contains artifacts and noise, and the low-contrast regions appear blurred. The images obtained with the proposed Deep Guess framework are more regular and have well-distinguished details. Finally, the last row of the figure displays the Deep Guess images computed by the neural network. We underline that these images are generated by state-of-the-art methods for CT reconstruction, and our proposal demonstrates that these methods can be successfully exploited to achieve higher-quality images with limited computational effort.

\begin{table}[]
    \centering
     \setlength\tabcolsep{0pt}
    \begin{tabular*}{\linewidth}{@{\extracolsep{\fill}} ll  cc ccc }
        \toprule
  &    \multicolumn{3}{c}{ $\bx^{(0)}$ } & \multicolumn{3}{c}{ Output } \\
        \cmidrule(lr){2-4} \cmidrule(lr){5-7} 
  &     &  SSIM & RE        &  SSIM & RE     & iters   \\

\midrule    
\multirowcell{5}{CP}
& zeros           & - & -                 & 91.72 & 0.0539 & 500      \\
& FBP             & 84.47 & 0.0744        & 91.08 & 0.0561 & 70      \\
& TV ($K=10$)     & 62.37 & 0.2535        & 82.70 & 0.0940 & 96  \\
& \textit{DG} by FBP-LPP    & 89.81 & 0.0560        & 91.87 & 0.0124 & 29   \\
& \textit{DG} by TV-RISING  & 81.76 & 0.1036        & 91.88 & 0.0606 & 66    \\
\midrule 
NETT & \quad -        & - & -                 & 90.25 & 0.0629 & 300     \\
\bottomrule
    \end{tabular*}
    \caption{Values of the quantitative metrics, computed on the Mayo Clinic test image reconstructions from the $\mathcal{G}_{360,360}$, with different initial guesses $\bx^{(0)}$ for the CP,  and with NETT method.}
    \label{tab:imm107mayo_G360360}
\end{table}

\begin{table}[]
    \centering
     \setlength\tabcolsep{0pt}
    \begin{tabular*}{\linewidth}{@{\extracolsep{\fill}} ll  cc ccc }
        \toprule
  &   \multicolumn{3}{c}{$\bx^{(0)}$ } & \multicolumn{3}{c}{ Output } \\
        \cmidrule(lr){2-4} \cmidrule(lr){5-7} 
  &     &  SSIM & RE        &  SSIM & RE     & iters   \\

\midrule
\multirowcell{5}{CP with \\ p=0.5}
& \quad zeros        & - & -               & 78.08 & 0.1039 & 500    \\
& \quad FBP         & 30.56 & 0.2955       & 77.01 & 0.1076 & 400\\
& \quad TV ($K=15$)   & 67.18 & 0.1906     & 75.36 & 0.1147 & 400\\
& \quad \textit{DG} by FBP-LPP    & 80.74 & 0.0876        & 83.16 & 0.0792 & 12    \\
& \quad \textit{DG} by TV-RISING  & 76.99 & 0.1059        & 78.79 & 0.0989 & 19    \\
\midrule
\multirowcell{5}{CP with \\ p=0.2 \\ }
& \quad zeros        & - & -               & 77.03 & 0.1071 & 500      \\
& \quad FBP       & 30.56 & 0.2955         & 76.35 & 0.1102 & 400     \\
& \quad TV ($K=15$)  & 67.18 & 0.1906      & 76.43 & 0.1093 & 400\\
& \quad \textit{DG} by FBP-LPP    & 80.74 & 0.0876        & 82.75 & 0.0794 & 18    \\
& \quad \textit{DG} by TV-RISING  & 76.99 & 0.1059        & 79.65 & 0.0958 & 36    \\
\midrule
NETT & \quad -        & - & -                 & 61.12 & 0.1576 & 300     \\
\bottomrule
    \end{tabular*}
    \caption{Values of the quantitative metrics, computed on the Mayo Clinic test image reconstructions for the $\mathcal{G}_{180,60}$ protocol, with different initial guesses $\bx^{(0)}$ for the CP and for $p=0.5$ and $p=0.2$,  and with NETT method.}
    \label{tab:imm107mayo_G18060}
\end{table}

\newpage
\subsection{Analysis of different Deep Guess implementations \label{subsec:coule}}

\begin{figure}

   
\footnotesize 
CP from  \hspace{25mm} CP from    \hspace{25mm} CP from   \hspace{24 mm} CP from  \\ 
zeros\hspace{30mm} DG by FBP-LPP     \hspace{15mm} DG by TV-RISING   \hspace{10mm} DG by \TpV-RISING  \\
    SSIM=92.73      \hspace{20mm} SSIM=98.67     \hspace{20mm} SSIM=98.05      \hspace{19mm} SSIM=98.04      \\
    \includegraphics[width=0.22\textwidth]{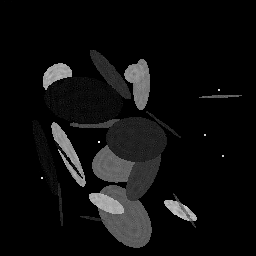}
    \includegraphics[width=0.22\textwidth]{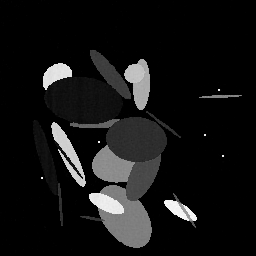}
    \includegraphics[width=0.22\textwidth]{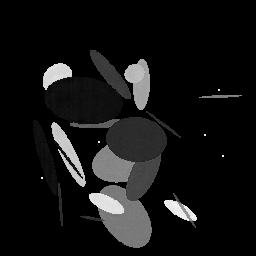}
    \includegraphics[width=0.22\textwidth]{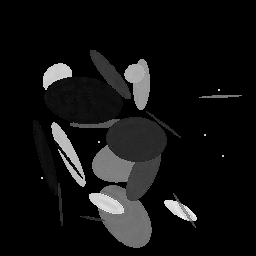}   \\ 
    \includegraphics[trim= 5mm 25mm 30mm 25mm, clip,width=0.22\textwidth]{ imm/coule14_zeros_IS_200iter.png}
    \includegraphics[trim= 5mm 25mm 30mm 25mm, clip,width=0.22\textwidth]{imm/coule14_FbpLpp_IS.png}
    \includegraphics[trim= 5mm 25mm 30mm 25mm, clip,width=0.22\textwidth]{imm/coule14_SgpRising_IS.png}
    \includegraphics[trim= 5mm 25mm 30mm 25mm, clip,width=0.22\textwidth]{imm/coule14_CpRising_IS.png}\\ 
    \includegraphics[trim= 22mm 5mm 13mm 45mm, clip,width=0.22\textwidth]{ imm/coule14_zeros_IS_200iter.png}
    \includegraphics[trim= 22mm 5mm 13mm 45mm, clip,width=0.22\textwidth]{imm/coule14_FbpLpp_IS.png}
    \includegraphics[trim= 22mm 5mm 13mm 45mm, clip,width=0.22\textwidth]{imm/coule14_SgpRising_IS.png}
    \includegraphics[trim= 22mm 5mm 13mm 45mm, clip,width=0.22\textwidth]{imm/coule14_CpRising_IS.png}  \\ \vspace{0mm}
\hspace{36mm} DG by FBP-LPP     \hspace{15mm} DG by TV-RISING   \hspace{10mm} DG by \TpV-RISING  \\ \vspace{0mm}
    \hspace{36mm} SSIM=93.97     \hspace{20mm} SSIM=92.57      \hspace{19mm} SSIM=75.80      \\ \vspace{0mm}
 \hspace{36mm} 
    \includegraphics[trim= 5mm 25mm 30mm 25mm,clip,width=0.22\textwidth]{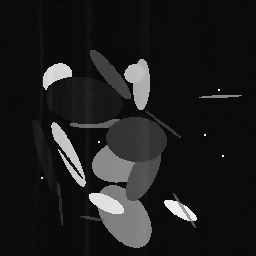}
    \includegraphics[trim= 5mm 25mm 30mm 25mm,clip,width=0.22\textwidth]{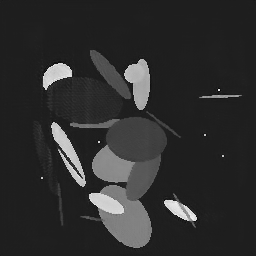}
    \includegraphics[trim= 5mm 25mm 30mm 25mm,clip,width=0.22\textwidth]{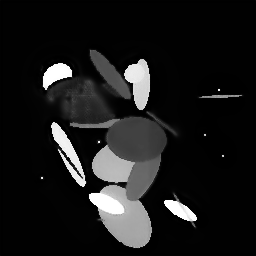}   \\ \vspace{0mm}
    \hspace{36mm} 
    \includegraphics[trim= 22mm 5mm 13mm 45mm, clip,width=0.22\textwidth]{imm/coule14_FbpLpp_DIG.png}
    \includegraphics[trim= 22mm 5mm 13mm 45mm, clip,width=0.22\textwidth]{imm/coule14_SgpRising_DIG.png}
    \includegraphics[trim= 22mm 5mm 13mm 45mm, clip,width=0.22\textwidth]{imm/coule14_CpRising_DIG.png} \\\vspace{1mm}
    
    \caption{  Results on the COULE test image shown in Figure \ref{fig:gt}, computed using the CP algorithm from the null image (first column) and some of the proposed Deep Guess acceleration approaches. The last row displays the zooms of the images generated by the Deep Guess neural networks within the proposed framework. The SSIM metric is reported for all images.}
    \label{fig:coule14}
\end{figure}

In this Section, we analyse and compare various possible implementations of the Deep Guess step on the COULE test set.
Here, we set $p=0.5$ and $K=10$ for both TV- and \TpV-based Deep Guesses. The regularization parameter $\lambda = 10^{-3}$ was heuristically determined to optimize performance when the CP solver is applied from the zeros $\bx^{(0)}$. This parameter remains unchanged for all the experiments of this section.\\
Table \ref{tab:medieCOULE} presents the results obtained on the COULE test set images. All SSIM and RE values are averages calculated across the test samples, with standard deviations provided in parentheses. The first two columns pertain to the Deep Guess image, the subsequent metrics correspond to the final reconstructions, and the last column indicates the range of iterations required for the CP method to converge in the second step of our scheme. \\
We observe that the ground truth-free RISING strategies for computing Deep Guess perform very well, losing, in the final output, only one point of SSIM compared to the corresponding LPP approaches (96.17 against 97.32 for the TV approach  and 93.01 against 94.54 for the \TpV\ one).Additionally, the TV-LPP and TV-RISING frameworks outperform their \TpV\ counterparts, offering higher metrics and, more importantly, smaller standard deviations, making them more reliable. Regarding the results limited to the Deep Guess step to compute $\bx^{(0)}$ (in the left part of the table), we note that \TpV-RISING performs significantly worse than the others. We attribute this behavior to the potential use of poor local minima as targets during network training.
\begin{table}[]
    \centering \footnotesize 
    \setlength\tabcolsep{0pt}
    \begin{tabular*}{\linewidth}{@{\extracolsep{\fill}} l  cc ccc }
    
    \toprule
    &   \multicolumn{2}{c}{ $\bx^{(0)}$ } & \multicolumn{3}{c}{ Output } \\
        \cmidrule(lr){2-3} \cmidrule(lr){4-6} 
       &  SSIM &  RE        & SSIM & RE     &  iters   \\
    \midrule    
zeros          & -$\quad$  & -$\quad$                              & 88.70 (3.37)  & 0.0851 (0.0073)  & 200      \\
\textit{DG} by FBP-LPP     &90.02 (4.19) & 0.0642 (0.0163)              & 96.20 (2.34)  & 0.0321 (0.0057)  & 45-55   \\
\textit{DG} by TV-LPP      & 87.45 (3.28) & 0.0527 (0.0074)              & 97.32 (1.03)  & 0.0290 (0.0025)  & 20-32    \\
\textit{DG} by TV-RISING   & 90.79 (2.60) & 0.0619 (0.0105)              & 96.17 (2.49)  & 0.0341 (0.0047)  & 47-55    \\
\textit{DG} by \TpV-LPP     & 93.80 (2.07) & 0.1372 (0.0517)              & 94.54 (4.23)  & 0.0530 (0.0120)  & 75-80     \\
\textit{DG} by \TpV-RISING  & 83.78 (10.84)& 0.1434 (0.0515)              & 93.01  (4.49)  & 0.0538 (0.0142)  & 73-85    \\
    \bottomrule
    \end{tabular*}

    \caption{Mean values (standard deviations) of the quantitative metrics, computed on the COULE test images. The first two columns present the metric values for the computed image $\bx^{(0)}$ (DG step); on the right, the metrics for the final reconstructed image are shown.}
    \label{tab:medieCOULE}
\end{table}
In Figure \ref{fig:coule14}, we present the reconstruction (together with the two considered crops) of the test image displayed in Figure \ref{fig:gt}. The considered  solutions are those computed from the zero-image  (first column) and from the Deep Guess computed by  the FBP-LPP (second column), the TV-RISING (third column) and the \TpV-RISING (last column) networks. The corresponding Deep Guess images are displayed at the bottom.    \\
The left image exhibits ringing artifacts, particularly in the homogeneous regions, as clearly shown in the two magnified zooms-in. In contrast, all the other final solutions computed from the DG images are clear and accurate. 
In the context of the Deep Guess $\bx^{(0)}$ images displayed in the final rows, the \TpV-RISING method exhibits the lowest accuracy, evidenced by blurred ellipses and the loss of low-contrast objects, aligning with the metrics reported in Table \ref{tab:medieCOULE}. Although the FBP-LPP image achieves the highest SSIM value, it displays streaking artifacts in low-contrast regions, a common issue associated with analytical reconstructions in sparse CT.

\begin{figure}
    \centering
    \includegraphics[width=0.9\textwidth]{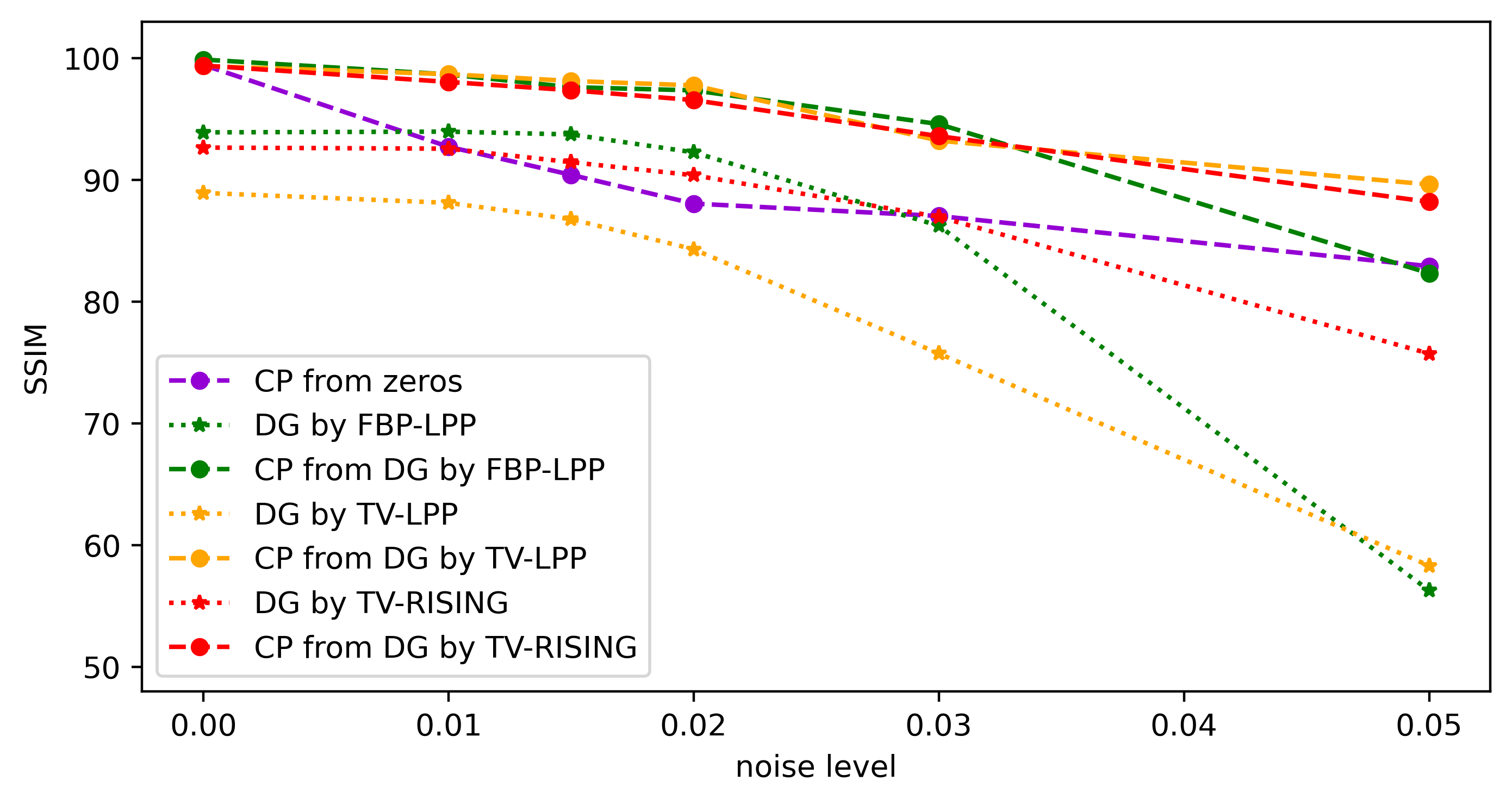}
    \caption{Plot of the SSIM values for the COULE test image reconstructions at different noise levels. Dotted lines with star markers represent the outputs of the Deep Guess neural networks, while dashed lines with circle markers correspond to the final images produced by the proposed framework and its competitors.}
    \label{fig:COULEStability}
\end{figure}

To assess the stability with respect to the noise of the hybrid methods and their corresponding CP solutions, we consider sinograms affected by noise at varying levels.
Now, we only focus on the FBP-LPP, TV-LPP, and the TV-RISING schemes.
In Figure \ref{fig:COULEStability}, we plot the SSIM values as a function of the noise level $\nu$ ranging from 0 (noiseless data) to 0.05. The results are relative to the previously considered test image, and we use the same parameter settings tuned for the previous experiments so that the tests shown in Figure \ref{fig:coule14} correspond to the case where the noise level is 0.01 in this plot. We plot the SSIM  both fro the DG image (dotted lines with star markers) and for the final reconstructions (dashed lines with circle markers).
As observed in \cite{evangelista2023rising}, the RISING network (red dotted line) is slightly less accurate, when tested on the same noise used for training, but more stable than the networks trained on the ground truth images (FBP-LPP and TV-LPP networks described by green and yellow dotted lines, respectively). In fact, when the noise level is similar to that analyzed during the training phase (i.e. $\nu=0.01$), the LPP  generally provides higher values. However, there is a turning point around the noise level of $\nu =0.03$, where the noise perturbation significantly impacts the data, and RISING outperforms LPP (as the red dotted line surpasses the yellow and orange ones). Regarding the final reconstructions, it is evident that the DG Accelerated framework consistently outperforms other methods across all the experiments considered (indicated by the red and orange dashed lines). Notably, even though the CP solver initialized from a zero image (represented by the purple dashed line) achieves SSIM values around 90 for low noise levels, its performance does not significantly degrade under conditions of high noise.

\section{Conclusion}\label{sec:conclusion}

This paper introduces a novel hybrid framework integrating a deep learning tool in a model-based regularized approach for reconstructing CT images from sparse views. 
The framework addresses a non-convex model that enforces sparsity in the image gradient domain and employs deep learning to achieve accurate image restoration with minimal computational effort.\\
The proposed approach is composed of two steps. The first one, called Deep Guess,  builds upon state-of-the-art efficient methods that generate output images using a post-processing convolutional neural network applied to a coarse reconstruction. To ensure data coherence and achieve the final reconstruction, in the second step we incorporate  a few iterations of a non-convex Model-Based Iterative Reconstruction  method.
We propose several implementations of Deep Guess, including those that do not rely on ground-truth data, thereby improving the practical applicability of our framework for real-world scenarios.\\
Experiments conducted on test problems using real CT images across various sparse geometries demonstrate the effectiveness of our approach. It surpasses other methods based on the same CP algorithm with different initial iterates, as well as the state-of-the-art hybrid methods (that merge variational and data-driven approaches) in terms of both accuracy and computational complexity, being 15 to 25 times faster. A thorough analysis of the stability of the proposed framework with respect to data noise demonstrates that the ground-truth-free implementations maintain strong generalization performance even as noise levels increase.\\
Additionally,  we emphasize the versatility of our framework, which can be tailored to specific tasks in both Deep Guess computation and the selection of MBIR methods for final image reconstruction, as depicted in Figure \ref{fig:graph_abstr}. 
Finally, our study demonstrates that deep learning tools can serve as a useful foundation for model-based approaches. By integrating these tools, we can accelerate computations while preserving the mathematical properties and characterizations of the final reconstructions.

\section{Acknowledgements and funding}
This work has been partially supported by Gruppo Nazionale per il Calcolo Scientifico (GNCS-INdAM) within the project ``Deep Variational Learning: un approccio combinato per la ricostruzione di immagini'', and by the PRIN 2022 project ``STILE: Sustainable Tomographic Imaging with Learning and rEgularization'', project  code: 20225STXSB, funded by the European Commission under the NextGeneration EU programme.

\bibliographystyle{abbrv}
\bibliography{biblio_cut}

\begin{thebibliography}{10}

\bibitem{antun2020instabilities}
V.~Antun, F.~Renna, C.~Poon, B.~Adcock, and A.~C. Hansen.
\newblock On instabilities of deep learning in image reconstruction and the
  potential costs of ai.
\newblock {\em Proceedings of the National Academy of Sciences},
  117(48):30088--30095, 2020.

\bibitem{baguer2020computed}
D.~O. Baguer, J.~Leuschner, and M.~Schmidt.
\newblock Computed tomography reconstruction using deep image prior and learned
  reconstruction methods.
\newblock {\em Inverse Problems}, 36(9):094004, 2020.

\bibitem{beck2009fast}
A.~Beck and M.~Teboulle.
\newblock A fast iterative shrinkage-thresholding algorithm for linear inverse
  problems.
\newblock {\em SIAM journal on imaging sciences}, 2(1):183--202, 2009.

\bibitem{beck2009fastwavelet}
A.~Beck and M.~Teboulle.
\newblock A fast iterative shrinkage-thresholding algorithm with application to
  wavelet-based image deblurring.
\newblock In {\em 2009 IEEE International Conference on Acoustics, Speech and
  Signal Processing}, pages 693--696. IEEE, 2009.

\bibitem{boink2019learned}
Y.~E. Boink and C.~Brune.
\newblock Learned svd: solving inverse problems via hybrid autoencoding.
\newblock {\em arXiv preprint arXiv:1912.10840}, 2019.

\bibitem{bonettini2008scaled}
S.~Bonettini, R.~Zanella, and L.~Zanni.
\newblock A scaled gradient projection method for constrained image deblurring.
\newblock {\em Inverse problems}, 25(1):015002, 2008.

\bibitem{borsdorf2008wavelet}
A.~Borsdorf, R.~Raupach, T.~Flohr, and J.~Hornegger.
\newblock Wavelet based noise reduction in ct-images using correlation
  analysis.
\newblock {\em IEEE transactions on medical imaging}, 27(12):1685--1703, 2008.

\bibitem{candes2006robust}
E.~J. Cand{\`e}s, J.~Romberg, and T.~Tao.
\newblock Robust uncertainty principles: exact signal reconstruction from
  highly incomplete frequency information.
\newblock {\em IEEE Transactions on information theory}, 2006.

\bibitem{candes2008enhancing}
E.~J. Cand\`es, M.~B. Wakin, and S.~P. Boyd.
\newblock Enhancing sparsity by reweighted {L}1 minimization.
\newblock {\em Journal of Fourier analysis and applications}, 2008.

\bibitem{cascarano2022plug}
P.~Cascarano, E.~Loli~Piccolomini, E.~Morotti, and A.~Sebastiani.
\newblock Plug-and-play gradient-based denoisers applied to {CT} image
  enhancement.
\newblock {\em Applied Mathematics and Computation}, 2022.

\bibitem{caselles2015total}
V.~Caselles, A.~Chambolle, and M.~Novaga.
\newblock Total variation in imaging.
\newblock {\em Handbook of mathematical methods in imaging}, 2015.

\bibitem{chambolle2011first}
A.~Chambolle and T.~Pock.
\newblock A first-order primal-dual algorithm for convex problems with
  applications to imaging.
\newblock {\em Journal of mathematical imaging and vision}, 2011.

\bibitem{chan2005recent}
T.~Chan et~al.
\newblock Recent developments in total variation image restoration.
\newblock {\em Mathematical Models of Computer Vision}, 2005.

\bibitem{chartrand2007exact}
R.~Chartrand.
\newblock Exact reconstruction of sparse signals via nonconvex minimization.
\newblock {\em IEEE Signal Processing Letters}, 2007.

\bibitem{chen2017low}
H.~Chen, Y.~Zhang, M.~K. Kalra, F.~Lin, Y.~Chen, P.~Liao, J.~Zhou, and G.~Wang.
\newblock Low-dose ct with a residual encoder-decoder convolutional neural
  network.
\newblock {\em IEEE transactions on medical imaging}, 36(12):2524--2535, 2017.

\bibitem{daubechies2010iteratively}
I.~Daubechies, R.~DeVore, M.~Fornasier, and C.~S. G{\"u}nt{\"u}rk.
\newblock Iteratively reweighted least squares minimization for sparse
  recovery.
\newblock {\em Communications on Pure and Applied Mathematics: A Journal Issued
  by the Courant Institute of Mathematical Sciences}, 63(1):1--38, 2010.

\bibitem{demirel2021p}
G.~Demirel, {\.I}.~Yildirim, and M.~Erta{\c{s}}.
\newblock P norm regularized breast tomosynthesis image reconstruction.
\newblock In {\em 2021 Medical Technologies Congress (TIPTEKNO)}, pages 1--4.
  IEEE, 2021.

\bibitem{donoho2006compressed}
D.~L. Donoho.
\newblock Compressed sensing.
\newblock {\em IEEE Transactions on information theory}, 2006.

\bibitem{evangelista2023rising}
D.~Evangelista, E.~Morotti, and E.~Loli~Piccolomini.
\newblock {RISING}: A new framework for model-based few-view {CT} image
  reconstruction with deep learning.
\newblock {\em Computerized Medical Imaging and Graphics}, 2023.

\bibitem{evangelista2023ambiguity}
D.~Evangelista, E.~Morotti, E.~L. Piccolomini, and J.~Nagy.
\newblock Ambiguity in solving imaging inverse problems with
  deep-learning-based operators.
\newblock {\em Journal of Imaging}, 9(7):133, 2023.

\bibitem{jin2017deep}
K.~H. Jin, M.~T. McCann, E.~Froustey, and M.~Unser.
\newblock Deep convolutional neural network for inverse problems in imaging.
\newblock {\em IEEE transactions on image processing}, 26(9):4509--4522, 2017.

\bibitem{li2020nett}
H.~Li, J.~Schwab, S.~Antholzer, and M.~Haltmeier.
\newblock Nett: Solving inverse problems with deep neural networks.
\newblock {\em Inverse Problems}, 36(6):065005, 2020.

\bibitem{piccolomini2018reconstruction}
E.~Loli~Piccolomini, V.~L. Coli, E.~Morotti, and L.~Zanni.
\newblock Reconstruction of 3d {X}-ray {CT} images from reduced sampling by a
  scaled gradient projection algorithm.
\newblock {\em Computational Optimization and Applications}, 2018.

\bibitem{piccolomini2016fast}
E.~Loli~Piccolomini and et~al.
\newblock A fast total variation-based iterative algorithm for digital breast
  tomosynthesis image reconstruction.
\newblock {\em Journal of Algorithms \& Computational Technology}, 2016.

\bibitem{loli2021model}
E.~Loli~Piccolomini and E.~Morotti.
\newblock A model-based optimization framework for iterative digital breast
  tomosynthesis image reconstruction.
\newblock {\em Journal of Imaging}, 2021.

\bibitem{lv2020nonlocal}
D.~Lv, Q.~Zhou, J.~K. Choi, J.~Li, and X.~Zhang.
\newblock Nonlocal tv-gaussian prior for bayesian inverse problems with
  applications to limited ct reconstruction.
\newblock {\em Inverse Problems \& Imaging}, 14(1), 2020.

\bibitem{mccollough2016tu}
C.~McCollough.
\newblock Tu-fg-207a-04: Overview of the low dose ct grand challenge.
\newblock {\em Medical physics}, 43(6Part35):3759--3760, 2016.

\bibitem{mehranian2013x}
A.~Mehranian, M.~R. Ay, A.~Rahmim, and H.~Zaidi.
\newblock X-ray ct metal artifact reduction using wavelet domain $ l\_
  $\{$0$\}$ $ sparse regularization.
\newblock {\em IEEE transactions on medical imaging}, 32(9):1707--1722, 2013.

\bibitem{monga2021algorithm}
V.~Monga, Y.~Li, and Y.~C. Eldar.
\newblock Algorithm unrolling: Interpretable, efficient deep learning for
  signal and image processing.
\newblock {\em IEEE Signal Processing Magazine}, 2021.

\bibitem{morotti2021green}
E.~Morotti and et~al.
\newblock A green prospective for learned post-processing in sparse-view
  tomographic reconstruction.
\newblock {\em Journal of Imaging}, 2021.

\bibitem{morotti2024space}
E.~Morotti, D.~Evangelista, A.~Sebastiani, and E.~L. Piccolomini.
\newblock Space-variant total variation boosted by learning techniques in
  few-view tomographic imaging.
\newblock {\em arXiv preprint arXiv:2404.16900}, 2024.

\bibitem{persson2001total}
M.~Persson, D.~Bone, and H.~Elmqvist.
\newblock Total variation norm for three-dimensional iterative reconstruction
  in limited view angle tomography.
\newblock {\em Physics in Medicine \& Biology}, 46(3):853, 2001.

\bibitem{porta2015new}
F.~Porta, M.~Prato, and L.~Zanni.
\newblock A new steplength selection for scaled gradient methods with
  application to image deblurring.
\newblock {\em Journal of Scientific Computing}, 65(3):895--919, 2015.

\bibitem{purisha2017controlled}
Z.~Purisha, J.~Rimpel{\"a}inen, T.~Bubba, and S.~Siltanen.
\newblock Controlled wavelet domain sparsity for x-ray tomography.
\newblock {\em Measurement Science and Technology}, 29(1):014002, 2017.

\bibitem{ritschl2011improved}
L.~Ritschl and et~al.
\newblock Improved total variation-based {CT} image reconstruction applied to
  clinical data.
\newblock {\em Physics in Medicine \& Biology}, 2011.

\bibitem{rodriguez2013total}
P.~Rodr{\'\i}guez.
\newblock Total variation regularization algorithms for images corrupted with
  different noise models: a review.
\newblock {\em Journal of Electrical and Computer Engineering}, 2013(1):217021,
  2013.

\bibitem{ronneberger_2015_unet}
O.~Ronneberger, P.~Fischer, and T.~Brox.
\newblock U-net: Convolutional networks for biomedical image segmentation.
\newblock In N.~Navab, J.~Hornegger, W.~M. Wells, and A.~F. Frangi, editors,
  {\em Medical Image Computing and Computer-Assisted Intervention -- MICCAI
  2015}, pages 234--241, Cham, 2015. Springer International Publishing.

\bibitem{rudin1992nonlinear}
L.~I. Rudin, S.~Osher, and E.~Fatemi.
\newblock Nonlinear total variation based noise removal algorithms.
\newblock {\em Physica D: nonlinear phenomena}, 60(1-4):259--268, 1992.

\bibitem{schnurr2019simulation}
A.-K. Schnurr, K.~Chung, T.~Russ, L.~R. Schad, and F.~G. Z{\"o}llner.
\newblock Simulation-based deep artifact correction with convolutional neural
  networks for limited angle artifacts.
\newblock {\em Zeitschrift f{\"u}r Medizinische Physik}, 29(2):150--161, 2019.

\bibitem{shu2022sparse}
Z.~Shu and A.~Entezari.
\newblock Sparse-view and limited-angle ct reconstruction with untrained
  networks and deep image prior.
\newblock {\em Computer Methods and Programs in Biomedicine}, 226:107167, 2022.

\bibitem{sidky2014cttpv}
E.~Y. Sidky and et~al.
\newblock Constrained ${T}p{V}$ minimization for enhanced exploitation of
  gradient sparsity: Application to {CT} image reconstruction.
\newblock {\em IEEE Journal of Translational Engineering in Health and
  Medicine}, 2014.

\bibitem{sidky2021docnn}
E.~Y. Sidky and et~al.
\newblock Do {CNN}s solve the {CT} inverse problem?
\newblock {\em IEEE Transactions on Biomedical Engineering}, 2021.

\bibitem{sidky2012convex}
E.~Y. Sidky, J.~H. J{\o}rgensen, and X.~Pan.
\newblock Convex optimization problem prototyping for image reconstruction in
  computed tomography with the {C}hambolle--{P}ock algorithm.
\newblock {\em Physics in Medicine \& Biology}, 2012.

\bibitem{sidky2008image}
E.~Y. Sidky and X.~Pan.
\newblock Image reconstruction in circular cone-beam computed tomography by
  constrained, total-variation minimization.
\newblock {\em Physics in Medicine \& Biology}, 53(17):4777, 2008.

\bibitem{tian2011low}
Z.~Tian and et~al.
\newblock Low-dose {CT} reconstruction via edge-preserving total variation
  regularization.
\newblock {\em Physics in Medicine \& Biology}, 2011.

\bibitem{tong2018edge}
S.~Tong, B.~Han, and J.~Tang.
\newblock Edge-guided tvp regularization for diffuse optical tomography based
  on radiative transport equation.
\newblock {\em Inverse Problems}, 34(11):115009, 2018.

\bibitem{van2016fast}
W.~Van~Aarle, W.~J. Palenstijn, J.~Cant, E.~Janssens, F.~Bleichrodt,
  A.~Dabravolski, J.~De~Beenhouwer, K.~Joost~Batenburg, and J.~Sijbers.
\newblock Fast and flexible x-ray tomography using the astra toolbox.
\newblock {\em Optics express}, 24(22):25129--25147, 2016.

\bibitem{van2015astra}
W.~Van~Aarle, W.~J. Palenstijn, J.~De~Beenhouwer, T.~Altantzis, S.~Bals, K.~J.
  Batenburg, and J.~Sijbers.
\newblock The astra toolbox: A platform for advanced algorithm development in
  electron tomography.
\newblock {\em Ultramicroscopy}, 157:35--47, 2015.

\bibitem{venkatakrishnan2013plug}
S.~V. Venkatakrishnan, C.~A. Bouman, and B.~Wohlberg.
\newblock Plug-and-play priors for model based reconstruction.
\newblock In {\em 2013 IEEE global conference on signal and information
  processing}, pages 945--948. IEEE, 2013.

\bibitem{vogel1996iterative}
C.~R. Vogel and M.~E. Oman.
\newblock Iterative methods for total variation denoising.
\newblock {\em SIAM Journal on Scientific Computing}, 1996.

\bibitem{wang2021non}
J.~Wang.
\newblock Non-convex lp regularization for sparse reconstruction of electrical
  impedance tomography.
\newblock {\em Inverse Problems in Science and Engineering}, 29(7):1032--1053,
  2021.

\bibitem{wang2018nonconvex}
S.~Wang, I.~Selesnick, G.~Cai, Y.~Feng, X.~Sui, and X.~Chen.
\newblock Nonconvex sparse regularization and convex optimization for bearing
  fault diagnosis.
\newblock {\em IEEE Transactions on Industrial Electronics}, 65(9):7332--7342,
  2018.

\bibitem{wang2003multiscale}
Z.~Wang, E.~P. Simoncelli, and A.~C. Bovik.
\newblock Multiscale structural similarity for image quality assessment.
\newblock In {\em The Thrity-Seventh Asilomar Conference on Signals, Systems \&
  Computers, 2003}, volume~2, pages 1398--1402. Ieee, 2003.

\bibitem{xiang2021fista}
J.~Xiang, Y.~Dong, and Y.~Yang.
\newblock Fista-net: Learning a fast iterative shrinkage thresholding network
  for inverse problems in imaging.
\newblock {\em IEEE Transactions on Medical Imaging}, 40(5):1329--1339, 2021.

\bibitem{xu2016accelerated}
Q.~Xu, D.~Yang, J.~Tan, A.~Sawatzky, and M.~A. Anastasio.
\newblock Accelerated fast iterative shrinkage thresholding algorithms for
  sparsity-regularized cone-beam ct image reconstruction.
\newblock {\em Medical physics}, 43(4):1849--1872, 2016.

\bibitem{superultra}
S.~Ye, Z.~Li, M.~T. McCann, Y.~Long, and S.~Ravishankar.
\newblock Unified supervised-unsupervised (super) learning for x-ray ct image
  reconstruction.
\newblock {\em IEEE Transactions on Medical Imaging}, 40(11):2986--3001, 2021.

\bibitem{zhu2018image}
B.~Zhu, J.~Z. Liu, S.~F. Cauley, B.~R. Rosen, and M.~S. Rosen.
\newblock Image reconstruction by domain-transform manifold learning.
\newblock {\em Nature}, 555(7697):487--492, 2018.

\end{thebibliography}
\end{document}